\documentclass[12pt]{amsart} 
\usepackage{amssymb}

\makeatletter
\let\over\@@over
\let\atop\@@atop
\let\above\@@above
\let\overwithdelims\@@overwithdelims
\let\atopwithdelims\@@atopwithdelims
\let\abovewithdelims\@@abovewithdelims

\def\cases#1{\left\{\,\vcenter{\normalbaselines\m@th
    \ialign{$##\hfil$&\quad##\hfil\crcr#1\crcr}}\right.}
\def\matrix#1{\null\,\vcenter{\normalbaselines\m@th
    \ialign{\hfil$##$\hfil&&\quad\hfil$##$\hfil\crcr
      \mathstrut\crcr\noalign{\kern-\baselineskip}
      #1\crcr\mathstrut\crcr\noalign{\kern-\baselineskip}}}\,}

\makeatother

\newtheorem{lemma}{Lemma}
\newtheorem{corollary}{Corollary}

\theoremstyle{definition}

\newtheorem{remarks}{Remarks}

\def\eqlabel#1{\label{eq#1}}
\def\eqref#1{(\ref{eq#1})}

\let\<\langle
\let\>\rangle

\makeatother

\begin{document}

\title{A sharp Schwarz inequality on the boundary}
\author{Robert Osserman}
\address{}
\curraddr{MSRI\\
1000 Centennial Drive\\
Berkeley, CA  94720-5070\\
USA}

\begin{abstract}
A number of classical results reflect the fact that if a holomorphic 
function maps the unit disk into itself taking the origin into the
origin,
 and if some boundary point $b$ maps to the boundary, then the
map is a magnification at $b$. We prove a sharp quantitative version of
this result which also sharpens a classical result of Loewner, and
which implies that the map is a strict magnification at 
$b$ unless it
is a rotation.
\end{abstract}

\thanks{Research at MSRI is supported in part by NSF grant DMS-9701755
}

\maketitle

The standard Schwarz Lemma states that an analytic function $f(z)$
mapping the unit disk into itself, with $f(0)=0$, must map each
smaller disk $|z|<r<1$ into itself and (as a result) satisfy
$|f'(0)|\leq 1$. Furthermore, unless $f$ is a rotation, one has strict
inequality $|f'(0)|<1$ and $f$ maps each disk $|z|\leq r<1$ into a 
strictly smaller one.

It is an elementary consequence of Schwarz' Lemma that if $f$ extends
continuously to some boundary point $b$ with $|b|=1$, and if $|f(b)|=1$
and $f'(b)$ exists, then $|f'(b)|\geq 1$. David Gilbarg raised the
question if one again has strict inequality unless $f$ is a rotation. The
answer is ``yes'' (Lemma 1 below), but it does not follow from the standard
Schwarz inequality. One needs a stronger form (Lemma 2 below) where one
has a quantitative bound on how \emph{much} each disk $|z|\leq r<1$ is shrunk
if $f$ is not a rotation.

\begin{lemma}[The boundary Schwarz Lemma]
Let $f(z)$ satisfy
\begin{itemize}
\item[(a)] $f(z)$ is analytic for $|z|<1$,
\item[(b)] $|f(z)|<1$ for $|z|<1$,
\item[(c)] $f(0)=0$,
\item[(d)] for some $b$ with $|b|=1$, $f(z)$ extends continuously to
$b$, $|f(b)|=1$, and $f'(b)$ exists. 
\end{itemize}

Then
\begin{equation}
\eqlabel{1}
|f'(b)|\geq \frac{2}{1+|f'(0)|}.
\end{equation} 
\end{lemma}

\begin{corollary}
Under hypotheses (a) - (d),
\begin{equation}
\eqlabel{2}
|f'(b)|\geq 1
\end{equation}
and
\begin{equation}
\eqlabel{3}
|f'(b)|>1 \mbox{ unless $f(z)=e^{i\alpha}z$, $\alpha$ real.}
\end{equation}
\end{corollary}

\begin{proof}
Inequalities (2) and (3) follow immediately from (1) together with 
the standard Schwarz Lemma.
\end{proof}

\begin{corollary}
Let $f$ satisfy conditions (a), (b), (c) of the lemma, and suppose
that f extends continuously to an arc $C$ on $|z|=1$, with $|f(z)|=1$
on $C$. Then the length $s$ of $C$ and the length $\sigma$ of $f(C)$
satisfy
\begin{equation}
\eqlabel{4}
\sigma\geq \frac{2}{1+|f'(0)|}\,
s.
\end{equation}
\end{corollary}

\begin{proof}
By the reflection principle, $f$ extends to be analytic on the interior
of $C$ and therefore satisfies condition (d) of Lemma 1. Hence (4)
follows from (1).
\end{proof}

\begin{remarks}
\begin{itemize}
\item[1.]
Again by the standard Schwarz Lemma, (4) implies that $\sigma \geq s$,
and $\sigma > s$ unless f is a rotation. That is the content of a
classical theorem of Loewner \cite{L}. (See also Velling \cite{V}.)
\item[2.]
The length $\sigma$
 of $f(C)$ is to be taken with multiplicity, if $f(C)$
is a multiple covering of the image.
\item[3.]
Inequality (1) is sharp, with equality possible for each value of $|f'(0)|$.
\item[4.]
One can drop the condition (c) that $f(0)=0$. Analogous results hold for
any value of $f(0)$. See Lemma 3, the General Boundary Lemma, below.
\item[5.]
One does not need to assume that $f$ extends continuously to $b$. For
example, if $f$ has a radial limit $c$ at $b$, with $|c|=1$, and if
$f$ has a radial derivative at $b$, then that derivative also satisfies
the inequality (1). More generally, if for some $b$ with $|b|=1$, there
exists a sequence $z_n$ such that $z_n \rightarrow b$ and 
$f(z_n) \rightarrow c$, with $|c|=1$, then
\begin{equation} 
\eqlabel{5}
\underline{\lim}_{z_n \rightarrow
b}\left|\frac{f(z_n)-c}{|z_n|-|b|}\right|
\geq \underline{\lim}_{z_n \rightarrow b}\,\frac{1-|f(z_n)|}{1-|z_n|}
\geq \frac{2}{1+|f'(0)|}.
\end{equation}
Both Lemma 1 and the statement about radial limits are immediate
consequences, since in either case we may choose $z_n=t_nb$ for 
$t_n$ real, $t_n \rightarrow 1$, and the left-hand side of (5)
becomes $|f'(b)|$.
\end{itemize}
\end{remarks}

\begin{lemma}[Interior Schwarz Lemma]
Let $f(z)$ satisfy conditions (a), (b), (c) of Lemma 1. Then
\begin{equation}
\eqlabel{6}
|f(z)| \leq |z| \frac{|z|+|f'(0)|}{1+|f'(0)||z|} \mbox{ for $|z|<1$}.
\end{equation}
\end{lemma}

\begin{proof}
Let $g(z)= \frac{f(z)}{z}$. Then by the standard Schwarz Lemma, either
$f$ is a rotation, or else $|g(z)|<1$ for $|z|<1$. In the former case,
$|f'(0)|=1$ and (6) holds trivially. So we need only consider the
second case, where $|g(z)|<1$. Furthermore, since inequality (6) is
unaffected by rotations, we may assume that $g(0)=f'(0)=a$, where
$0 \leq a <1$. Then (6) is equivalent to
\begin{equation}
\eqlabel{7}
|g(z)| \leq \frac{|z|+a}{1+a|z|} \mbox{ for $|z|<1$,  with $a=g(0)$}.
\end{equation}
But that is an immediate consequence of the standard Schwarz-Pick
version of the Schwarz Lemma, which says that $g$ must map each
disk $|z|<r$ into the image of that disk under the linear fractional
map 
$$
G(z)= \frac{z+a}{1+az}
$$
which is a circular disk whose diameter is the interval
$$
\Bigl[\frac{a-r}{1-ar}, \frac{a+r}{1+ar}\Bigl]
$$
of the real axis.

Hence,
$$
|z|=r \Rightarrow |g(z)| \leq \frac{a+r}{1+ar} = \frac{|z|+a}{1+a|z|},
$$
which proves (7), and hence (6).
\end{proof}

\begin{remarks}
\begin{itemize}
\item[1.]
For related sharpened forms of the interior Schwarz Lemma, see 
Mercer \cite{M}.
\item[2.]
Inequality (7) is sharp, with equality for $g(z)=G(z)$, $z=r$.
Hence, inequality (6) is sharp, with equality for the function
$$
f(z)=z \frac{z+a}{1+az},\   0 \leq a<1,
$$
when  $z$ is on the positive real axis. The same function gives
equality in (1) when $b=1$.
\item[3.]
When $f$ is not a rotation, (6) is a strict improvement on the
standard Schwarz Lemma, since the second factor on the right is
strictly less than 1 when $|f'(0)|<1$.
\end{itemize}
\end{remarks}

\begin{proof}[Proof of Lemma 1.]
Let f satisfy conditions (a),(b),(c) of Lemma 1. Then, using the
upper bound (6) for $|f(z)|$, we have for any $b$ and $c$ with
$|b|=1$, $|c|=1$, 
$$
\left|\frac{f(z)-c}{|z|-|b|}\right| \geq \frac{1-|f(z)|}{1-|z|} \geq 
\frac{1+|z|}{1+|f'(0)||z|}.
$$
As $|
z|\rightarrow 1$, the right-hand side tends to $\frac{2}{1+|f'(0)|}.$

This proves (5), and as noted in Remark 5 above, Lemma 1 follows.
\end{proof}

\begin{lemma}[The General Boundary Lemma]
\label{3}
Under hypotheses (a),(b), and (d) of Lemma 1, one has
\begin{equation}
\eqlabel{8}
|f'(b)| \geq \frac{2}{1+|F'(0)|} \frac{1-|f(0)|}{1+|f(0)|},
\end{equation}
where $F$ is defined in $(9)$
below and satisfies $|F'(0)| \leq 1$,
with strict inequality unless $F$ is a rotation and $f$ is an automorphism
of the unit disk. 
\end{lemma}
\begin{proof}
Let
\begin{equation}
\eqlabel{9}
F(z)= \frac{f(z)-f(0)}{1- \overline{f(0)}f(z)}.
\end{equation}
Then $F$ satisfies the hypotheses of Lemma 1, and therefore
\begin{equation}
\eqlabel{10}
|F'(b)| \geq \frac{2}{1+|F'(0)|}
,
\end{equation}
with equality if and only $F$ is a rotation. But a calculation gives
$$
F'(z)=f'(z) \frac{1-|f(0)|^{2}}{[1- \overline{f(0)} f(z)]^{2}}.
$$
Since $|f(b)|=1$ implies
$$
|1- \overline{f(0)} f(b)| \geq 1-|\overline{f(0)}f(b)| = 1-|f(0)|,
$$
we have 
\begin{equation}
\eqlabel{11}
|F'(b)|=|f'(b)| \frac{1-|f(0)|^{2}}{|1-\overline{f(0)}f(b)|^{2}}
\leq|f'(b)|\frac{1+|f(0)|}{1-|f(0)|}.
\end{equation}
Combining (10) and (11) yields (8).
\end{proof}

\begin{remarks}[Concluding Remarks]
\begin{itemize}
\item[1.]
An interesting special case of Lemma 1 is when $f'(0)=0$, in which
case inequality (1) implies $|f'(b)|\geq 2$. Clearly equality holds
for
\begin{equation}
\eqlabel{12}
f(z)=e^{i\alpha}z^2, \mbox{\quad $\alpha$ real}.
\end{equation}
Furthermore, that is the only case of equality; the same type of
argument used to prove Lemmas 1 and 2 yields a stronger inequality
that implies $|f'(b)|>2$ unless $f$ is of the form (12). More generally, 
the argument of the standard Schwarz lemma shows that if
 $f(z)=\sum_{n=0}^{\infty} a_n z^n$ satisfies (a), (b) of Lemma 1
and if
\begin{equation}
\eqlabel{13}
a_0=a_1= \cdots a_{k-1}=0,
\end{equation}
then $|a_k|\leq 1$, and $|a_k|=1$ if and only if
\begin{equation}
\eqlabel{14}
f(z)= e^{i\alpha}z^k , \mbox{\quad $\alpha$ real}.
\end{equation}
Furthermore, either (14) holds, or else $|f(z)|<|z|^k$ for $|z|<1$.
The argument of Lemma 2 yields the stronger result that
\begin{equation}
\eqlabel{15}
|f(z)| \leq |z|^k \frac{|z|+|a_k|}{1+|a_k||z|}.
\end{equation}
Using (15) in the proof of Lemma 1 then shows that if also condition (d)
of Lemma 1 holds, then
\begin{equation}
\eqlabel{16}
|f'(b)| \geq k + \frac{1-|a_k|}{1+|a_k|}.
\end{equation}
It follows that $|f'(b)| \geq k$, with equality only if $f$ is of the
form (14).
\item[2.]
A corollary of Lemma 3 is that under the same hypotheses, one has
\begin{equation}
\eqlabel{17}
|f'(b)| \geq \frac{1-|f(0)|}{1+|f(0)|}
\end{equation}
and the inequality is strict unless $f$ is an automorphism of the
unit disk. In this context, see Cara\-th\'eo\-dory \cite{C1}, pp.54-55,
on Julia's Theorem.
\item[3.]
For related results, and other types of boundary
Schwarz Lemmas, see Cara\-th\'eo\-dory \cite{C2}, pp.54-55, on the converse
of Julia's
Theorem, Pommerenke \cite{P}, p. 71, on the Julia-Wolff Lemma, and
the paper of Burns and Krantz \cite{BK}.
\end{itemize}
\end{remarks}


\begin{thebibliography}{}

\bibitem
[BK]{BK} Burns, D.M. and Krantz, S.G. ``Rigidity of holomorphic
mappings and a new Schwarz Lemma at the boundary.'' {\em J. Amer. Math.
Soc.} {\bf 7} (1994) 661-676. 

\bibitem[C1]{C1} Carath\'eodory, C. {\em Conformal Representation}
Cambridge University Press 1952.

\bibitem[C2]{C2} Carath\'eodory, C. {\em Theory of Functions} Vol.2,
New York, Chelsea 1954.

\bibitem[L]{L} L\"{o}wner. K. ``Untersuchungen \"{u}ber schlichte
konforme Abbildungen des Einheitskreises. I.'' {\em Math. Ann.}
{\bf 89} (1923) 103-121.

\bibitem[M]{M} Mercer, P.R. ``Sharpened versions of the Schwarz Lemma.''
{\em J. Math. Anal. Appl.} {\bf 205} (1997) 508-511.

\bibitem[P]{P} Pommerenke, Ch. {\em Boundary Behavior of Conformal Maps},
New York, Springer-Verlag 1992.

\bibitem[V]{V} Velling, J.A. ``The uniformizations of rectangles, an
exercise in Schwarz's Lemma.'' {\em Amer. Math. Monthly} {\bf 39} (1992)
112-115.



\end{thebibliography}
\end{document}